\documentclass[12pt,thmsa]{article}
\usepackage{amsfonts}

\input{tcilatex}
\begin{document}

\author{Karl-Georg Schlesinger \qquad \\
Erwin Schr\"{o}dinger Institute for Mathematical Physics\\
Boltzmanngasse 9\\
A-1090 Vienna, Austria\\
e-mail: kgschles@esi.ac.at}
\title{On a quantum analog of the Grothendieck-Teichm\"{u}ller group }
\date{}
\maketitle

\begin{abstract}
We introduce a self-dual, noncommutative, and noncocommutative Hopf algebra $%
\mathcal{H}_{GT}$ which takes for certain Hopf categories (and therefore
braided monoidal bicategories) a similar role as the
Grothendieck-Teichm\"{u}ller group for quasitensor categories. We also give
a result which highly restricts the possibility for similar structures for
higher weak $n$-categories ($n\geq 3$) by showing that these structures
would not allow for any nontrivial deformations. Finally, give an explicit
description of the elements of $\mathcal{H}_{GT}$.
\end{abstract}

\section{The Hopf algebra $\mathcal{H}_{GT}$}

In \cite{Dri} Drinfeld introduced the Grothendieck-Teichm\"{u}ller group by
considering the (formal) reparametrizations of the data (commutativity and
associativity isomorphisms) of a quasitensor category. Consider now braided
(weak) monoidal bicategories arising from the representations of a Hopf
category (as defined in \cite{CF}) on 2-vector spaces (see \cite{KV}), i.e.
on certain module categories. Let us assume, in addition, that the Hopf
category itself is given as the category of finite dimensional
representations of a quasi-trialgebra, satisfying a quasitriangularity and
coquasitriangularity condition. This is analogous to understanding the
Grothendieck-Teichm\"{u}ller group $GT$ as a universal symmetry of
quasitriangular quasi-Hopf algebras (see e.g. \cite{CP}) which via their
category of finite dimensional representations then give rise to the afore
mentioned quasitensor categories. Considering the question of a universal
symmetry of quasitriangular quasi-Hopf algebras is e.g. of interest in the
study of symmetries of moduli spaces of two dimensional conformal field
theories (see \cite{Kon}) since two dimensional conformal field theories are
closely linked to the heighest weight representation of quasitrinagular
quasi-Hopf algebras through their vertex algebras. Since the author has
shown that any three dimensional extended topological quantum field theory
in the sense of \cite{KL} uniquely determines a trialgebra (\cite{Sch}) and
these three dimensional extended topological quantum field theories are
supposed to be related to two dimensional boundary conformal field theories,
the question of a universal symmetry of trialgebras is of potential interest
to the question of symmetries on moduli spaces of two dimensional boundary
conformal field theories.

\bigskip

\begin{remark}
Note that the above mentioned result, linking trialgebras to extended
topological quantum field theories, also shows that the restriction of the
consideration to Hopf categories which are representation categories of a
trialgebra still includes a large and - form the physics persepective - the
most important class of examples of such structures.
\end{remark}

\bigskip

Let us begin by commenting on some of the involved notions or give
references to the relevant literature, respectively. Especially, we will
introduce the notion of a trialgebra in detail, now. For the notion of
(quasi-) Hopf algebras, quasitriangularity and coquasitriangularity, etc.,
we refer to any of the many excellent introductions to Hopf algebras and
quantum groups, available now (e.g. \cite{CP} or \cite{KS}). The notion of
quasitensor category which is used in Drinfeld's definition of $GT$ is given
as a category together with a tensor product $\otimes $ on it where $\otimes 
$ need not be symmetric but satisfying a commutativity constraint ``up to
isomorphism''. For the purpose of this article, the reader should imagine a
quasitensor category simply as the category of finite dimensional
representations of a quasitriangular quasi-Hopf algebra and the
commutativity constraint to be given by a universal $R$-matrix. For the a
detailed introduction of quasitensor categories and their link to Hopf
algebras and quantum groups, we refer to \cite{CP}.

Let us next introduce the concept of a trialgebra:

\bigskip

\begin{definition}
A trialgebra $(A,*,\Delta ,\cdot )$ with $*$ and $\cdot $ associative
products on a vector space $A$ (where $*$ may be partially defined, only)
and $\Delta $ a coassociative coproduct on $A$ is given if both $(A,*,\Delta
)$ and $(A,\cdot ,\Delta )$ are bialgebras and the following compatibility
condition between the products is satisfied for arbitrary elements $%
a,b,c,d\in A$: 
\[
(a*b)\cdot (c*d)=(a\cdot c)*(b\cdot d)
\]
whenever both sides are defined.
\end{definition}

\bigskip

Trialgebras were first suggested in \cite{CF} as an algebraic means for the
construction of four dimensional topological quantum field theories. It was
observed there that the representation categories of trialgebras have the
structure of so called Hopf algebra categories (see \cite{CF}) and it was
later shown explicitly in \cite{CKS} that from the data of a Hopf category
one can, indeed, construct a four dimensional topological quantum field
theory. The first explicit examples of trialgebras were constructed in \cite
{GS1} and \cite{GS2} by applying deformation theory, once again, to the
function algebra on the Manin plane and some of the classical examples of
quantum algebras and function algebras on quantum groups. In \cite{GS4} it
was shown that one of the trialgebras constructed in this way appears as a
symmetry of a two dimensional spin system. Besides this, the same trialgebra
can also be found as a symmetry of a certain system of infinitely many
coupled $q$-deformed harmonic oscillators.

\bigskip

\begin{definition}
We call a trialgebra quasitriangular (coquasitriangular) if one of the
bialgebras contained in it is quasitriangular (coquasitriangular). We call a
trialgebra $\left( A,\cdot ,*,\Delta \right) $ biquasitriangular if $\left(
A,\cdot ,\Delta \right) $ is quasitriangular and $\left( A,*,\Delta \right) $
is coquasitriangular and if for the $R$-matrix $R_{\cdot }$ of $\left(
A,\cdot ,\Delta \right) $ and the linear form $R_{*}$ expressing the
coquasitriangularity of $\left( A,*,\Delta \right) $, the following
condition holds: 
\[
\left[ R_{\cdot },\widehat{R_{*}}\right] =0
\]
where $\widehat{R_{*}}$ is the $R$-matrix of a Hopf algebra dual to $\left(
A,*,\Delta \right) $.
\end{definition}

\bigskip

We will speak of a quasi-trialgebra if there is a Drinfeld coassociator $%
\alpha $ for $\Delta $, and one of the two products has a dual associator $%
\beta $ such that $\alpha $ and $\beta $ satisfy a similar commutator
condition as the $R$-matrices above (we will study this condition in detail
in the next section).

\bigskip

\begin{lemma}
The (formal) reparametrizations of the data of the above mentioned
biquasitriangular quasi-trialgebras define a self-dual noncommutative and
noncocommutative Hopf algebra $\mathcal{H}_{GT}$.
\end{lemma}

\proof%
In a quasitriangular quasi-Hopf algebra we have two kinds of data which are
transformed by $GT$ as a universal symmetry (see \cite{Dri}): The $R$-matrix
and the coassociator $\alpha $. In the precise definition of $GT$ the
completion of the transformations of these data with respect to a cartain
class of formal power series is considered (see \cite{Dri} or \cite{CP} for
a comprehensive introduction to $GT$). In a biquasitriangular
quasi-trialgebra we have four types of data: The two matrices $R_{\cdot }$
and $R_{*}$, the coassociator $\alpha $, and the associator $\beta $. We ask
for the universal symmetry given by transformations of these data (including
the same completion with respect to formal power series as in the case of $%
GT $), now.

First, observe that on the data $\left( R_{\cdot },\alpha \right) $ taken
alone $GT$ acts just by definition. Considering formal linear combinations
of the data $\left( R_{\cdot },\alpha \right) $, we can, obviously, extend
this to an action of the group algebra of $GT$ (which naturally has the
structure of a Hopf algebra, see e.g. \cite{CP} or \cite{KS}). Second, the
class of all data $\left( R_{*},\beta \right) $ is dual to the class of all
data $\left( R_{\cdot },\alpha \right) $. So, concerning a universal
symmetry of the data $\left( R_{*},\beta \right) $ taken alone, we have to
have a dual of the action of $GT$, again. By the definition of the data $%
\left( R_{*},\beta \right) $, we can not have an action of a group there but
have to describe a universal symmetry by a coaction of a Hopf algebra. By
the above argument, this has to be the function algebra on $GT$ (where we
define the appropriate function algebra as the algebra of polynomial
functions, since $GT$ is a projective limit of algebraic groups and the
explicit definition of $GT$ in \cite{Dri} assures that the product of $GT$
correctly transforms into a coproduct as one proves by calculation from the
defining relations).

In consequence, if we would transform the data $\left( R_{\cdot },\alpha
\right) $ and $\left( R_{*},\beta \right) $ of the two bialgebras included
in a trialgebra separately, forgetting about the compatibility condition for
the two products of a trialgebra, the universal symmetry would be described
by the Drinfeld double $\mathcal{D}\left( GT\right) $ of $GT$, i.e. the
tensor product of the group algebra of $GT$ and the algebra of functions on $%
GT$. In the next step, we have to restrict to those transformations of the
complete set of data $R_{\cdot },R_{*},\alpha ,\beta $ which transform a
biquasitriangular quasi-trialgebra into a biquasitriangular
quasi-trialgebra. Obviously, this is a subspace of $\mathcal{D}\left(
GT\right) $. One proves by calculation from the compatibility relation of
the two products that it is a sub-Hopf algebra $\mathcal{H}_{GT}$, indeed.

It remains to show that $\mathcal{H}_{GT}$ is self-dual, noncommutative, and
noncocommutative: The self-duality follows from the fact that the classes of
data $\left( R_{\cdot },\alpha \right) $ and $\left( R_{*},\beta \right) $
are dual to each other. $\mathcal{H}_{GT}$ can not be commutative since one
of the factors of $\mathcal{D}\left( GT\right) $ restricts to the group
algebra of $GT$ and $GT$ is non-abelian. Finally, $\mathcal{H}_{GT}$ is
noncocommutative, then, since it is self-dual. This completes the proof. 
\endproof%

\bigskip

\begin{lemma}
There is an algebra morphism from $\mathcal{H}_{GT}$ to the group algebra of
the Grothendieck-Teichm\"{u}ller group.
\end{lemma}

\proof%
Since, as shown above, $\mathcal{H}_{GT}$ is a sub-Hopf algebra of $\mathcal{%
D}\left( GT\right) $ and one of the factors of $\mathcal{D}\left( GT\right) $
is just the group algebra of $GT$, the conclusion follows. 
\endproof%

\bigskip

\begin{remark}
Observe that the above map is only a morphism with respect to the
associative alegbra structure of $\mathcal{H}_{GT}$. Also, it can not be
surjective since there are compatibility constraints between the
commutativity and associativity isomorphisms and their dual structures.
\end{remark}

\bigskip

One could have the idea to extend this approach to higher braided weak
monoidal weak $n$-categories beyond the level of bicategories where for
tricategories one would expect an algebraic structure in the form of a
vector space equipped with two associative products and two coassociative
coproducts to generate these tricategories via representation theory. We
will call such an algebraic structure a quadraalgebra. The compatibilities
are given by requiring that any of the coproducts together with the two
products defines a trialgebra plus the requirement that the two coproducts
are compatible by the dual relation to the compatibility relation for the
two products.

\bigskip

\begin{lemma}
There do not exist nontrivial deformations of a trialgebra into a
quadraalgebra.
\end{lemma}

\proof%
With similar arguments as given above one can show that the universal
symmetry of a quasi-quadraalgebra, with two quasitriangularity conditions
and one coquasitriangularity condition satisfied, is given by a trialgebra $%
\mathcal{T}_{GT}$. where there is an algebra morphism from $\mathcal{T}_{GT}$
- as an associative algebra - to the group algebra of $GT$. Besides this,
one proves that both bialgebras (which are Hopf algebras, even, in this
case) within $\mathcal{T}_{GT}$ are self-dual and the two products have to
agree (this follows, again from the symmetry of the classes of data on which 
$\mathcal{T}_{GT}$ acts as a universal symmetry). Both products are, as a
consequence, universally defined, then. Since both products agree, we can
without loss of generality assume that we have a unital product. But then we
can apply an Eckmann-Hilton type argument to conclude that the product is
abelian. So, $\mathcal{T}_{GT}$ is a commutative and cocommutative self-dual
Hopf algebra. Besides this, $\mathcal{T}_{GT}$ is a sub-Hopf algebra of $%
\mathcal{D}\left( GT\right) $. But because of the algebra morphism from $%
\mathcal{T}_{GT}$ to the group algebra of $GT$, given by the previuos lemma, 
$\mathcal{T}_{GT}$ is determined by an abelian subgroup of $GT$, then. But
by definition of $GT$ (see \cite{Dri}), we get triviality, then, i.e. $%
\mathcal{T}_{GT}$ consists - up to rational factors - of the identity, only.
But triviality of $\mathcal{T}_{GT}$ means that we can not have a nontrivial
formal deformation theory of quasi-quadraalgebras with suitable
(co)quasitriangularity conditions. Remenbering that in the definition of $GT$
the coassociator is the essential part of the data (see \cite{Dri} where
this is already noted), we can extend this conclusion to general
quasi-quadraalgebras. This concludes the proof. 
\endproof%

\bigskip

So, on the level of tricategories arising via representation theory from
quadraalgebras (and, consequently, for higher categorical levels linked to
corresponding higher algebras), one gets only generalizations of braided
monoidal structures which do not allow for deformations. Especially, as we
have just shown, there is no nontrivial deformation theory of trialgebras
into quadraalgebras, further generalizing the deformation of groups into
Hopf algebras into trialgebras. So, on the level of trialgebras a kind of
stability is reached. Observe that this non existence of deformations is
much stronger than usual rigidity in cohomology theory since we can not only
exclude deformations in a given category of structures but also deformations
to higher categorical analogs of the structure. E.g. the usual rigidity
results in the theory of classical Lie algebras do not exclude the
deformation of the universal envelope into a noncommutative and
noncocommutative Hopf algebra but - as we just mentioned - we can exclude
deformations of trialgebras into algebraic structures involving four or more
products and coproducts joined in a compatible way. We suggest the term 
\textit{ultrarigidity} for this kind of stability.

\bigskip

\begin{remark}
Since Hopf categories are linked to four dimensional topological field
theory (as Hopf algebras are to the three dimensional case), see \cite{CF}
and \cite{CKS}, this seems on the algebraic level to mirror the fact that
geometry in dimension five and higher is in some sense much simpler than the
three and four dimensional cases. There is also a more physical
interpretation of this result: Since bialgebra categories are linked
supposedly to certain types of quantum field theories on noncommutative
spaces (see \cite{GS3}), we can see this as saying that quantum field theory
on noncommutative spaces is a stable structure in some sense, not allowing
for a further generalization of the passage from classical to quantum field
theory to quantum field theory on noncommutative spaces.
\end{remark}

\bigskip

So far, we have seen only a few abstract properties of the Hopf algebra $%
\mathcal{H}_{GT}$. We will give a more explicit description of $\mathcal{H}%
_{GT}$ in the next section.

\bigskip

\section{The explicit structure of $\mathcal{H}_{GT}$}

In \cite{Dri} an explicit description of the Grothendieck-Teichm\"{u}ller
group $GT$ is derived from the general definition of the group of
transformations of the associator and the braiding of a quasitensor
category. In this section, we want to do the same for the Hopf algebra $%
\mathcal{H}_{GT}$.

Recall that the elements of $GT$ can be written in the form $\left( \lambda
,f\right) $ with $\lambda \in \Bbb{Q}$ and $f$ belongs to the $\Bbb{Q}$%
-pro-unipotent completion of the free group of two generators where the
pairs $\left( \lambda ,f\right) $ satisfy certain conditions (see \cite{Dri}%
). Remember also that $f$ arises in the following way from the general
definition of $GT$: If we change the associator 
\[
\left( U\otimes V\right) \otimes W\rightarrow U\otimes \left( V\otimes
W\right) 
\]
this means multiplying it by an automorphism of $\left( U\otimes V\right)
\otimes W$. It can be shown that any such automorphism is of the form 
\[
f\left( \sigma _1^2,\sigma _2^2\right) \left( \sigma _1\sigma _2\right)
^{3n} 
\]
with $n\in \Bbb{Z}$ and $f$ as above. Here, $\sigma _1$, $\sigma _2$ are the
generators of the braid group $B_3$.

Now, assume that we have a Hopf category (see \cite{CF}) with associativity
isomorphism $\alpha $ for the tensor product and coassociativity isomorphism 
$\beta $ for the functorial coproduct. Note that while the possible
transformations of $\alpha $ are represented by automorphisms of tensor
products $\left( U\otimes V\right) \otimes W$, the possible
cotransformations of $\beta $ are of a dual nature and can formally be seen
as elements of the algebraic dual of the underlying vector space of the
automorphism group 
\[
Aut\left( \left( U\otimes V\right) \otimes W\right) 
\]
(remember that a Hopf category is, especially, $\Bbb{C}$-linear). So,
cotransformations of $\beta $ can be - up to linear combinations - written
in the form $\widehat{g}$ where $g$ is the second component of an element of 
$GT$ and\ $\widehat{}\ $ denotes the dualization operation as defined above.

\bigskip

Next, remember that the structure of $GT$ is basically determined by the
transformations of the associator (see \cite{Dri}, \cite{Kon}), i.e. in the
sequel we will forget about the component $\lambda $ coming from the
braiding.

In conclusion, we can describe the Hopf algebra $\mathcal{H}_{GT}$ as a
sub-Hopf algebra of the tensor product of the function algebra of $GT$ with
the Hopf algebra dual of $GT$ (as defined above), i.e. as a sub-Hopf algebra
of the Drinfeld double of $GT$.

In order to determine the concrete nature of this sub-Hopf algebra, we have
to use the compatibility condition between $\alpha $ and $\beta $ involved
in the definition of a Hopf category. While for objects $U,V,W$ $\alpha $ is
represented as an isomorphism 
\[
\varphi :\left( U\otimes V\right) \otimes W\rightarrow U\otimes \left(
V\otimes W\right) 
\]
$\beta $ is, again, given by an element $\widehat{\psi }$ of the dual of the
space of such transformations. The natural compatibility condition is, then, 
\begin{equation}
\psi ^{-1}\varphi =\varphi ^{-1}\psi  \label{B1}
\end{equation}
In order to assure that a pair $\left( f,\widehat{g}\right) $ transforms a
Hopf category into a Hopf category, we have to require that for the
transformed isomorphisms the above equation holds, too. Since $f$ and $g$
act on $\varphi $ and $\psi $, respectively, by the multiplication 
\begin{eqnarray*}
\varphi &\mapsto &\varphi f \\
\psi &\mapsto &\psi g
\end{eqnarray*}
it follows that 
\begin{equation}
g^{-1}\psi ^{-1}\varphi f=f^{-1}\varphi ^{-1}\psi g  \label{B2}
\end{equation}
Let 
\[
\chi =\psi ^{-1}\varphi 
\]
i.e. equation (\ref{B2}) reads as 
\begin{equation}
g^{-1}\chi f=f^{-1}\chi ^{-1}g  \label{B3}
\end{equation}
We are searching for a universal structure (i.e. not dependent on the choice
of Hopf category) of $\mathcal{H}_{GT}$, so, we have to require that (\ref
{B3}) holds for all possible choices of $\chi $.

Since equation (\ref{B1}) implies that 
\[
\chi ^2=1 
\]
i.e. $\chi $ is a projector, it follows that equation (\ref{B3}) holds for
all possible choices of $\chi $ iff it holds for $\chi $ being the identity.
So, equation (\ref{B2}) is equivalent to the condition 
\begin{equation}
g^{-1}f=f^{-1}g  \label{B4}
\end{equation}
for the elements $f,g$ of $GT$, i.e. we have a universal condition
determining $\mathcal{H}_{GT}$.

\bigskip

\begin{remark}
The symmetry inherent in the condition (\ref{B4}) is, of course, the source
of the self-duality of $\mathcal{H}_{GT}$.
\end{remark}

\begin{remark}
Mixed Tate motives over $Spec\left( \Bbb{Z}\right) $ are believed to be
given as representations of $GT$ (see \cite{Kon}). The explicit nature of
condition (\ref{B4}), in principle, allows for explicit calculations of the
quantum analogs of such motives as representations of $\mathcal{H}_{GT}$
which are also corepresentations of $\mathcal{H}_{GT}$. Given pairs of
representations of $GT$, one can use (\ref{B4}) to determine such
representations of $\mathcal{H}_{GT}$. On the other hand, the condition also
shows that one has to expect that representations of pairs $\left( f,%
\widehat{g}\right) $ satisfying (\ref{B4}) exist where neither the component 
$f$, nor $g$, derives from a full representation of $GT$, i.e. one has to
expect quantum motives which do not derive from a classical counterpart.
E.g. partial representations of $GT$ which would develop singularities, if
one would try to extend them to a full one, could play a role, here.
\end{remark}

\bigskip

We want to conclude this section with another small observation: In \cite{KL}
an algebraic framework - so called extended topological quantum field
theories - is developed in detail which allows for the inclusion of the case
of boundary conformal field theories into the algebraic description. It is
shown there that such theories are determined by modular categories $%
\mathcal{C}$ (i.e. certain quasitensor categories) together with a Hopf
algebra object $H$ \textit{in} $\mathcal{C}$. One can immediately define a
category Rep$\left( H\right) $ of representations of $H$ \textit{in} $%
\mathcal{C}$ from this.

\medskip

\begin{lemma}
The possible compatible transformations of Rep$\left( H\right) $ together
with $\mathcal{C}$ are determined by pairs $\left( f,g\right) $ of elements $%
f,g$ of $GT$ satisfying condition (\ref{B4}).
\end{lemma}

\proof%
Direct consequence of the definition of Rep$\left( H\right) $.%
\endproof%

\bigskip

So, from the view of the abstract quantum symmetry $\mathcal{H}_{GT}$, the
algebraic formulation of boundary conformal field theories given by \cite{KL}
and the structure of trialgebras and Hopf categories are just different
concrete realizations of one and the same quantum symmetry.

\bigskip

\begin{remark}
One can dually also formulate $\mathcal{H}_{GT}$ by starting from the Ihara
algebra \textit{Ih} (see \cite{Dri}, \cite{Iha1}, \cite{Iha2} for the
definition) instead of $GT$ (the Ihara algebra is closely related to the Lie
algebra of $GT$). The condition (\ref{B4}) translates then to the condition 
\begin{equation}
\left[ f,h\right] =0  \label{B5}
\end{equation}
for elements $f,h\in Ih$.

Remembering that the Lie algebra structure of \textit{Ih} derives - by
evaluation of the elements of \textit{Ih} on finite-dimensional metrized
(i.e. endowed with an invariant inner product) Lie algebra $g$ - from the
Kirillov bracket (see \cite{Dri} and for the definition of the Kirillov
bracket \cite{Kir}), condition (\ref{B5}) translates after evaluation on $g$
to 
\[
\left\{ f_g,h_g\right\} =0
\]
i.e. we can view it as requiring $h_g$ to behave as a symmetry relative to $%
f_g$ and vice versa.
\end{remark}

\bigskip

\section{Conclusion}

We have introduced a noncommutative analog $\mathcal{H}_{GT}$ of the
Grothendieck-Teichm\"{u}ller group in the form of a self-dual,
noncommutative, and noncocommutative Hopf algebra. Besides this, we have
given an explicit description of the elements of $\mathcal{H}_{GT}$. We also
proved a stability property (\textit{ultrarigidity}) excluding deformations
of certain higher categorical structures than bicategories. Further work
will deal, in particular, with physical applications of this stability
result.

\bigskip

\textbf{Acknowledgements:}

I thank H. Grosse for discussions on the topics involved and A. Goncharov
for very helpful explanations on classical motivic structures. Besides this,
I thank the Deutsche Forschungsgemeinschaft (DFG) for support by a research
grant and the Erwin Schr\"{o}dinger Institute for Mathematical Physics,
Vienna, for hospitality.

\end{document}